\documentclass[12pt,reqno]{amsart}

\usepackage{amsthm,amsmath,amssymb,amsfonts,amscd,latexsym,xypic}

\def\suffix{ps}
\newcount\system
\global\system=3   

\def\ifundefined#1{\expandafter\ifx\csname#1\endcsname\relax}
\ifundefined{figdir}\def\figdir{}\fi
\newcount\firstline
\newdimen\pswidth  \newdimen\xleft
\newdimen\psheight \newdimen\ytop \newdimen\ybot
\newcount\justx \newcount\justy
\global\justx=0 \global\justy=0
\newdimen\vpos \newtoks\labeL 
\newread\labeLfile \newdimen\xcoord \newdimen\ycoord
\newif\ifdoit 
\newbox\labox
\newdimen\xdvikwid 
\newdimen\xdvikht
\newdimen\pspoints
\newdimen\rwi
\newcount\temp
\def\readdim#1{\global\read\labeLfile to \temp
\global #1=\temp pt}
\def\figcrop#1{\par
\openin\labeLfile=\figdir#1.lbl                                              
\global\read\labeLfile to\firstline\message{#1}               
\global\read\labeLfile to\temp
\readdim{\ybot}
\readdim{\xleft}
\readdim{\ytop}
\global\read\labeLfile to\justx
\global\read\labeLfile to\justy
\global\read\labeLfile to\labeL
\readdim{\pswidth}
\global\advance\pswidth by -\xleft
\readdim{\psheight}
\global\advance\ybot by -\psheight
\global\advance\psheight by -\ytop
\global\read\labeLfile to\justx
\global\read\labeLfile to\justy
\global\read\labeLfile to\labeL
\vbox to\psheight{\vfill
\ifnum\system=1
\fi 
\ifnum\system=2
\fi 
\ifnum\system=3
                                                 \fi         
\ifnum\system=4
\fi         
\ifnum\system=1
\hbox to \pswidth{\kern-\xleft\special{postscriptfile \figdir#1.\suffix }\hfil}\fi
\ifnum\system=2
\hbox to \pswidth{\kern-\xleft\special{ps: plotfile \figdir#1.\suffix }\hfil}\fi
\ifnum\system=3
\hbox to \pswidth{\kern-\xleft\includegraphics{\figdir#1.\suffix}\hfil}\fi
\ifnum\system=4
\hbox to \pswidth{\kern-\xleft\includegraphics{\figdir#1.\suffix}\hfil}\fi
\ifnum\system=5
\hbox to \pswidth{\kern-\xleft\includegraphics{\figdir#1.\suffix}\hfil}\fi 
\ifnum\system=6
   \xdvikwid=\pswidth
   \xdvikht=\psheight
   {\global\divide\xdvikwid by \pspoints}
   {\global\divide\xdvikht by \pspoints}
   \rwi=\xdvikwid
    {\global\multiply\rwi by 10}
\hbox to \pswidth{\kern-\xleft\includegraphics{\figdir#1.\suffix\space}\hfil}\fi                   
\vskip -\baselineskip
\vskip -\ybot 
\vskip-\psheight %
\hbox to\pswidth  {\hss}%
\parindent=0pt\offinterlineskip                                       
\vpos=0 pt%
\loop\readdim{\xcoord}                                 
\ifdim \xcoord < -999pt \doitfalse\else\doittrue\fi                        
\ifdoit \advance \xcoord by -\xleft
\readdim{\ycoord}
\advance \ycoord by -\ytop                              
\global\read\labeLfile to\justx                                       
\global\read\labeLfile to\justy                                       
\global\read\labeLfile to\labeL
\global\setbox\labox=\hbox{\labeL\hskip-0.3em}%
\advance\vpos by-\ycoord                                              
\vskip-\vpos \vpos=\ycoord                                         
\hbox to\pswidth{\hskip\xcoord %
\hbox to 0pt{\ifnum\justx>0\hss\fi%
\vbox to0pt{%
\ifnum\justy<2\vss\fi%
\copy\labox\kern0pt%
\ifnum\justy>0\vss\fi}%
\ifnum\justx<2\hss\fi}%
\hss}%
\repeat%
\advance\vpos by-\psheight%
\vskip-\vpos %
}\closein\labeLfile}
\def\figplace#1#2#3{
\openin\labeLfile=\figdir#1.lbl
\ifeof \labeLfile
       \immediate\write16{***Can't find \figdir#1.lbl; Skipping it.***}
\else  \closein\labeLfile
       \null\hskip#2\raise #3 \hbox{\figcrop{#1}}
\fi
}
\def\bbbone{{\mathchoice {\rm 1\mskip-4mu l} {\rm 1\mskip-4mu l}
{\rm 1\mskip-4.5mu l} {\rm 1\mskip-5mu l}}}

\def\CC{\mathbb{C}}
\def\NN{\mathbb{N}}

\newcommand{\be}  {\begin{equation}}
\newcommand{\ee}  {\end{equation}}
\newcommand{\bea} {\begin{eqnarray}}
\newcommand{\eea} {\end{eqnarray}}
\newcommand{\lp}  {\left(}
\newcommand{\rp}  {\right)}

\newcommand{\cF}{\mathcal F}
\newcommand{\cJ}{\mathcal J}
\newcommand{\cC}{\mathcal C}
\newcommand{\cN}{\mathcal N}
\newcommand{\cA}{\mathcal A}
\newcommand{\cZ}{\mathcal Z}
\newcommand{\Om}{\Omega}
\newcommand{\ep}{\epsilon}
\newcommand{\si}{\sigma}
\newcommand{\ga}{\gamma}
\newcommand{\Ga}{\Gamma}
\newcommand{\al}{\alpha}
\newcommand{\de}{\delta}
\newcommand{\De}{\Delta}
\newcommand{\ph}{\phi}

\newcommand{\til}  {\tilde}
\def\Br{\overline}

\newcommand{\eqdef} {\stackrel{\rm def}{=}}

\newtheorem{Theorem}{Theorem}[section]
\newtheorem{Lemma}[Theorem]{Lemma}
\newtheorem{Proposition}[Theorem]{Proposition}
\newtheorem{Corollary}[Theorem]{Corollary}
\newtheorem{Remark}[Theorem]{Remark}

\newtheorem{Definition}[Theorem]{Definition}

\newcommand{\Z}{\mathbf Z}
\renewcommand{\P}{\mathbb P}        
\renewcommand{\O}{\mathcal O} 
\newcommand{\I}{\mathcal I}

\newcommand{\R}{\mathcal R} 
 
\newcommand{\F}{\mathcal F} 
\newcommand{\ux}{\mathbf x} 
\newcommand{\uy}{\mathbf y} 
 
\newcommand{\ra}{\rightarrow}

\newcommand{\lra}{\longrightarrow}

\renewcommand{\ker}{\text{ker}\,}

\newcommand{\demo}{\noindent {\sc Proof.}\;}

\begin{document}
\title{A regularity result for a locus of Brill type}
\author{Abdelmalek Abdesselam and Jaydeep Chipalkatti} 
\maketitle 

\parbox{12cm}{\small 
{\sc Abstract.} 
Let $n,d$ be a positive integers, with $d$ even (say $d=2e$). Write 
$N = \binom{n+d}{d}-1$, and let $X_{(n,d)} \subseteq \P^N$ 
denote the locus of degree $d$ hypersurfaces 
in $\P^n$ which consist of two $e$-fold hyperplanes. We calculate a 
bound on the Castelnuovo-regularity of its defining ideal, 
moreover we show that this variety is $r$-normal for $r \ge 2$. 
The latter part is proved by reducing the question to a 
combinatorial calculation involving Feynman diagrams and hypergeometric series. 
As such, it is a result of a tripartite collaboration of algebraic geometry, 
classical invariant theory, and modern theoretical physics.} 

\vspace{5mm} 

\mbox{\small AMS subject classification (2000): 05A15, 14F17, 14L35, 81T18.}

\medskip 

\parbox{12cm} 

\bigskip \bigskip 

\section{Introduction} 
The set of hypersurfaces of degree $d$ in $\P^n$ is parametrized by the 
projective space $\P^N$, where $N = \binom{n+d}{d}-1$. Assume that $d$ is 
even (say $d=2e$), and consider the subset of hypersurfaces 
which consist of two (possibly coincident) $e$-fold hyperplanes. 
In algebraic terms, this is the set of $(n+1)$-ary degree $d$ forms $F$
which factor as $F = L_1^e L_2^e$ for some linear forms $L_i$. 
This forms a subvariety of $\P^N$, which we will denote by $X_{(n,d)}$. 
It may be called of Brill type, in analogy with the following 
problem considered (and solved) by Brill: 
Find necessary and sufficient conditions for $F$ to factor 
into a product of linear forms, $F = L_1L_2 \dots L_d$. 

Let us write $X$ for $X_{(n,d)}$ if no confusion is likely. 
(We exclude the trivial case $n=1,d=2$ throughout.) 
Define  \[ m_0 = \lceil 2n+1 - \frac{n}{e} \rceil. \] 
Our main result is the following: 

\smallskip 

\begin{Theorem} \sl 
\begin{enumerate} 
\item[(i)]
The ideal $I_X$ is $m_0$-regular. 
A fortiori, it is minimally generated in degrees $\le m_0$. 
\item[(ii)]
The linear series on $X$ cut out by 
degree $r$ hypersurfaces in $\P^N$ is complete for $r \ge 2$.
\end{enumerate} 
\label{main.theorem} \end{Theorem}

\medskip 

The natural action of the group $SL_{n+1}$ on the 
imbedding $X \subseteq \P^N$ will be essential to the proof. 
As a byproduct we will get a formula for the image of 
each graded piece $(I_X)_r$ in the Grothendieck ring 
of $SL_{n+1}$-modules. 

Part (ii) of the theorem is the more delicate one. 
We reduce it to a question about transvectants of binary forms, and 
then resolve the issue using an explicit computation with Feynman diagrams. 
We hope that this technique should find a wider application. 

In section \S \ref{examples} we give classical invariant theoretic 
descriptions of the generators of $I_X$ for the cases 
$(n,d)=(1,8),(2,4)$. 

\begin{Remark} \rm It will be apparent that this is an instance of a 
problem which can be formulated rather generally. Given any partition 
$\lambda = (\lambda_1,\lambda_2,\dots)$ of $d$, one can define 
a subvariety of forms which factor as $\prod L_i^{\lambda_i}$. It is 
a natural problem to find $SL_{n+1}$-invariant equations for this variety; 
it is not completely settled even for binary forms. 

The case $\lambda=(d)$ corresponds to the Veronese imbedding (see \cite{JoeH}), and 
$\lambda=(1^d)$ to the Chow variety of degree $d$ zero cycles on $\P^n$ 
(see \cite[Ch.~4]{GKZ}). 
The case $\lambda=(e,e)$ is perhaps the next in order of complexity. A result for the 
case $(\lambda_1,\lambda_2)$ (with $\lambda_1 \neq \lambda_2)$ is under preparation. 
\end{Remark} 

\subsection{Preliminaries} 
The base field will be $\CC$. 
Let $V$ denote a complex vector space of dimension $n+1$, and 
write $W=V^*$. All subsequent constructions will be $SL(V)$-equivariant;
see \cite[Ch.~6 and 15]{FH} for the relevant representation theory. 
We will abbreviate $\text{Sym}^d \, V, \text{Sym}^r(\text{Sym}^d \, V)$ as 
$S_d,S_r(S_d)$ etc. 
If $\lambda$ is a partition, then $S_\lambda(-)$ will denote the associated 
Schur functor. All terminology from algebraic geometry follows \cite{Ha}. 

Fix a positive integer $d=2e$, and let $N = \binom{n+d}{d}-1$. 
Given the symmetric algebra 
\[ R = \bigoplus_{r \ge 0} S_r(S_d \, V), \] 
the space of degree $d$ hypersurfaces in $\P V$ is identified 
with 
\[ \P^N = \P \, S_d W = \text{Proj} \; R. \] 
Now define 
\begin{equation}
X_{(n,d)} = \{ [F] \in \P^N: F = (L_1 L_2)^e \;\; \text{for 
some $L_1,L_2 \in W$} \}. 
\end{equation}
This is an irreducible $2n$-dimensional projective subvariety of 
$\P^N$. 

\smallskip 

Recall the definition of regularity according to Mumford 
\cite[Ch.~6]{Mum}. 
\begin{Definition} \sl 
Let $\F$ be a coherent $\O_{\P^N}$-module, and $m$ an integer. 
Then $\F$ is said to be 
$m$-regular if $H^q(\P^N,\F(m-q)) = 0$ for $q \ge 1$. 
\end{Definition} 
It is known that $m$-regularity implies $m'$-regularity for 
all $m' \ge m$. 
Let $M$ be a graded $R$-module containing no submodules of 
finite length. Then (for the present purpose) we will say that $M$ is 
$m$-regular if its sheafification ${\widetilde M}$ is. In our case, 
$M=I_X$ (the saturated ideal of $X$), and ${\widetilde I_X} = \I_X$. 

We have the usual short exact sequence 
\begin{equation} 
0 \ra \I_X \ra \O_{\P^N} \ra \O_X \ra 0. 
\label{ses1} \end{equation} 
The map 
\begin{equation} \P W \times \P W \stackrel{f}{\lra} \P S_d \, W, \quad 
(L_1,L_2) \lra (L_1L_2)^e
\end{equation}
induces a natural isomorphism of $X$ with the quotient 
$(\P W \times \P W)// \Z_2$, and of the structure sheaf $\O_X$ with 
$(f_* \O_{\P W \times \P W})^{\Z_2}$. 

Using the Leray spectral sequence and the K{\"u}nneth formula, 
\[ \begin{aligned} 
{} & H^q(\P^N, f_* \O_{\P W \times \P W}(r)) = \\ 
\bigoplus\limits_{i+j=q} & 
H^i(\P W, \O_{\P^n}(re)) \otimes 
H^j(\P W, \O_{\P^n}(re)). 
\end{aligned} \] 
This group can be nonzero only in two cases: $i,j$ are either both 
$0$ or both $n$ (see \cite[Ch.~III,\S 5]{Ha}). 

\begin{Corollary} \sl 
We have an isomorphism $H^0(\O_X(r)) = S_2(S_{re})$ for 
$r \ge 0$. Moreover $H^{2n}(\O_X(r)) =0$ for $re \ge -n$. 
\end{Corollary} 

\section{The Proof of Theorem \ref{main.theorem}} 
Define the predicate 
\[ \R(q): \; H^q(\P^N, \I_X(m_0-q)) = 0. \] 
We would like to show $\R(q)$ for $q \ge 1$. 
Tensor the short exact sequence (\ref{ses1}) by $\O_{\P^N}(m_0-q)$ 
and consider the long exact sequence in cohomology. 
If $q \neq 1,2n+1$, then $\R(q)$ is immediate. By the choice of $m_0$, we have 
\[ e(m_0-2n-1) \ge -n. \] 
Hence $H^{2n}(\O_X(m_0-2n-1))=0$, which implies $\R(2n+1)$. 
Now $\R(1)$ is the case $r=m_0-1$ of the following result 
(which is part (ii) of the main theorem). 
\begin{Proposition} \sl 
Let $r \ge 2$. Then the morphism 
\[ \alpha_r: H^0(\O_{\P^N}(r)) \lra H^0(\O_X(r)) \] 
is surjective. 
\end{Proposition}
\demo 
The map $f$ can be factored as 
\[ \P W \times \P W \lra 
\P S_e \, W \times \P S_e \, W \lra \P S_d \, W. \] 
Tracing this backwards, we see that $\alpha_r$ is the composite 
\begin{equation}
S_r(S_d) \stackrel{1}{\lra} S_r(S_e \otimes S_e) 
\stackrel{2}{\lra} S_r(S_e) \otimes S_r(S_e) 
\stackrel{3}{\lra} S_{re} \otimes S_{re} 
\stackrel{4}{\lra} S_2(S_{re}), 
\label{alpha.r} \end{equation}
where 1 is given by applying $S_r(-)$ to the coproduct map, 2 is the  
projection coming from the `Cauchy decomposition' (see \cite{ABW}), 3 is the 
multiplication map, and 4 is the symmetrisation. 
Now we have a plethysm decomposition 
\begin{equation} H^0(\O_X(r)) = S_2(S_{re}) = 
\bigoplus\limits_p S_{(rd-2p,2p)}, 
\label{oxr} \end{equation}
where the direct sum is quantified over $0 \le p \le 
\lfloor \frac{re}{2} \rfloor$. Let $\pi_p$ denote the 
projection onto the $p$-th summand. Since any finite dimensional 
$SL(V)$-module is completely reducible, the cokernel of $\alpha_r$ is a 
direct summand of $H^0(\O_X(r))$. We will show that 
$\pi_p \circ \alpha_r \neq 0$ for any $p$, then Schur's lemma will 
imply that the cokernel is zero. 

The entire construction is functorial in $V$, hence if 
$U \subseteq V$ is any subspace, then the diagram 
\[ \diagram
S_r(S_d \, U) \dto \rto & S_{(rd-2p,2p)} \, U \dto \\
S_r(S_d \, V) \rto & S_{(rd-2p,2p)} \, V 
\enddiagram \] 
is commutative. If we further assume that $\dim U \ge 2$, then 
both vertical maps are injective. (Recall that $S_\lambda(V)$ vanishes
if and only if the number of parts in 
$\lambda$ exceeds $\dim V$.) Hence we may as well assume that 
$\dim V=2$. Thus we are reduced to the following statement: 

\begin{Proposition} \sl 
Assume $\dim V = 2$. Then 
the morphism $\pi_p \circ \alpha_r$ is nonzero for any 
$0 \le p \le \lfloor \frac{re}{2} \rfloor$. 
\label{prop.nonzero} \end{Proposition} 

The proof will be given in the sections \ref{transvectants} and \ref{Feynman}. 
Tentatively we will take the main theorem as proved, and 
interpose some examples. The following is a simple corollary to the theorem. 

\begin{Corollary} \sl 
In the Grothendieck ring of finite-dimensional $SL(V)$-modules, we 
have the equality
\[ [(I_X)_r] = 
[S_r(S_d)]- \sum\limits_{0 \le p \le \lfloor \frac{re}{2} \rfloor} 
[S_{(rd-2p,2p)}] \] 
Here $[-]$ denotes the formal character of a representation. 
\label{gr} \end{Corollary} 
\demo 
This follows because $(I_X)_r = \ker \alpha_r$.  \qed 

\smallskip 

Decomposing the plethysm $S_r(S_d)$ is in general a difficult 
problem. Explicit formulae are known only in very 
special cases -- see \cite{CGR,MacDonald} and the references therein. 
In particular the decomposition of $S_3(S_d)$ is given by Thrall's 
formula (see \cite{Plunkett1}), and then $(I_X)_3$ can be 
calculated in any specific case. 

\begin{Remark} \rm 
Note that $(I_X)_2 =0$, i.e., the ideal has no quadratic generators. 
If $n=1$, then $I_X$ is generated in degree $3$ and has 
regularity $3$, so its minimal resolution is linear. 
\end{Remark} 

\section{Examples} \label{examples} 
\subsection{Binary octavics} 
We will write down a complete set of invariant theoretic conditions 
necessary and sufficient for a degree eight binary form to lie in $X_{(1,8)}$. 
By what we have proved, the ideal has all of its generators in degree $3$. 
By Corollary \ref{gr} and Thrall's formula, 
\begin{equation} [(I_X)_3] = [ S_{18} \oplus S_{14} \oplus S_{12} 
\oplus S_{10} \oplus S_8 \oplus S_6]. 
\label{ix3} \end{equation}
Now (for instance) $S_{12}$ corresponds to a covariant of binary octavics 
of degree $3$ and order $12$. (This formalism is explained in \cite{ego4}). The 
module $S_3(S_8)$ contains two copies of $S_{12}$, so there is a two dimensional 
space of such covariants. An inspection reveals that 
\[ A = (F^2,F)_6, \quad B = ((F,F)_2,F)_4 \] 
can be taken as a basis for this space. (The fact that $A$ and $B$ 
have the right degree and order is clear from their definitions, so 
we only need to show that they are linearly independent for 
general $F$. This can be done by specializing $F$ to $x_0^6x_1^2 + x_0x_1^7$ 
and calculating directly.) 
Hence the required covariant must be $c_1 \, A + c_2 \, B$ for some 
constants $c_i$. Now specialize $F$ to 
$x_0^4x_1^4$, then by hypothesis the covariant vanishes. This gives a system 
of linear equations for the $c_i$, it has the solution 
$c_1:c_2 = 13:-63$. This determines the covariant (of course up to a scalar). 
By the same procedure, we can identify all the summands in 
(\ref{ix3}) as follows:
\begin{equation} \begin{array}{lll} 
(F^2,F)_3, & (F^2,F)_5, & 13 \, (F^2,F)_6 - 63 \, (F^2,F)_4, \\ 
(F^2,F)_7, & ((F,F)_6,F)_3, & 195 \, (F^2,F)_8 - 2744 \, ((F,F)_2,F)_6. 
\end{array} \label{oct} \end{equation}
We conclude that a form $F \in \P^8$ belongs to $X_{(1,8)}$, iff all the covariants in 
(\ref{oct}) vanish on $F$. 

\begin{Remark} \rm 
If $n \ge 2$, then $X_{(n,2)}$ is the variety of quadrics of 
rank at most $2$. It is a symmetric determinantal variety in the sense of 
\cite{JPW}, and its entire minimal resolution is deduced there. 
The ideal is generated in degree $3$ by the piece 
\[ (I_X)_3 = S_{(2,2,2)} \subseteq S_3(S_2). \] 
\end{Remark} 

\subsection{Ternary quartics} \label{ternary}
Assume $n=2,d=4$. By the main theorem, we know that the generators of 
$I_X$ lie in degrees $\le 4$. We will find them using an 
elimination theoretic computation. Define 
\[ \begin{aligned} 
   L_1 & = a_0 \, x_0 + a_1 \, x_1 + a_2 \, x_2, \quad 
   L_2 = b_0 \, x_0 + b_1 \, x_1 + b_2 \, x_2, \\
   F & = c_0 \, x_0^4 + c_1 \, x_0^3x_1 + \dots + c_{14} \, x_2^4;
\end{aligned} \] 
where the $a,b,c$ are indeterminates. 
Write $F = (L_1L_2)^2$ and then equate the coefficients of the monomials 
in $x_0,x_1,x_2$. This expresses 
each $c_i$ as a function of $a_0,\dots,b_2$, and hence defines a ring map 
\[ \CC \, [c_0,\dots,c_{14}] \lra \CC \, [a_0, \dots,b_2]. \] 
The kernel of this map is $I_X$. We calculated this in Macaulay-2, it 
turned out that in fact all the minimal generators are in degree $3$. 
By Corollary \ref{gr} and Thrall's formula, 
\[ (I_X)_3 = S_{(9,3)} \oplus S_{(6,0)} \oplus S_{(6,3)} 
\oplus S_{(4,2)} \oplus S_{(0,0)}. \] 
Now each summand corresponds to a concomitant of ternary quartics,
e.g., $S_{(9,3)}$ corresponds to one of degree $3$, order $6$ and 
class $3$. It is not difficult to identify the concomitants symbolically 
(see \cite{ego5} for the procedure), they are 
\begin{equation} \begin{array}{ll}
\alpha_x^2 \, \beta_x^3 \, \gamma_x \, (\alpha \, \gamma \, u)^2 
(\beta \, \gamma \, u), & 
\alpha_x^2 \, \beta_x^2 \, \gamma_x^2 \, (\alpha \, \beta \, \gamma)^2, \\ 
\alpha_x^2 \, \beta_x \, (\beta \, \gamma \, u)^2 
(\alpha \, \gamma \, u) (\alpha \, \beta \, \gamma), & 
\alpha_x \, \beta_x \, (\alpha \, \gamma \, u) (\beta \, \gamma \, u) 
(\alpha \, \beta \, \gamma)^2, \\ 
(\alpha \, \beta \, \gamma)^4. 
\end{array} \label{x2} \end{equation}
This is a rephrasing of the calculation in geometric terms: 
\begin{Theorem}\sl 
Let $F$ be a ternary quartic with zero scheme $C \subseteq \P^2$. 
Then $C$ consists of two (possibly coincident) double lines 
iff all the concomitants in (\ref{x2}) vanish on $F$. 
\end{Theorem}

\section{Transvectants} \label{transvectants} 
In this section we will break down Proposition \ref{prop.nonzero} into two 
separate questions about transvectants of binary forms. A general account of 
transvectants may be found in \cite{GrYo} and \cite{Olver}. 

We begin by describing the map $\alpha_r$ from (\ref{alpha.r}) in coordinates. 
(It is as yet unnecessary to assume $\dim V =2$.) 
Let 
\[ \ux^{(i)} = (x_0^{(i)}, \dots, x_n^{(i)}), \quad 
1 \le i \le r, \] 
be $r$ sets of $n+1$ variables, with their `copies' 
\[ \uy^{(i)} = (y_0^{(i)}, \dots, y_n^{(i)}), \quad 
1 \le i \le r. \] 

Let $F_i(\ux^{(i)}), 1 \le i \le r$ be degree $d$ forms, then the image 
$\alpha_r(\bigotimes\limits_{i=1}^r F_i)$ is calculated as follows: 
\begin{itemize} 
\item 
For each $F_i$, apply the polarization operator 
\[ \sum\limits_{\ell=0}^n y^{(i)}_\ell \frac{\partial}{\partial x^{(i)}_\ell} \] 
altogether $e$ times, denote the result by 
$F_i(\ux^{(i)},\uy^{(i)})$. 
\item 
Take the product $\prod\limits_i F_i(\ux^{(i)},\uy^{(i)})$, and make 
substitutions 
\[ x^{(i)}_\ell = x_\ell, \quad y^{(i)}_\ell = y_\ell, \] 
for all $i,\ell$. (This is tantamount to `erasing' the upper indices.)
This gives a form having degree $re$ each in $\ux,\uy$, which is 
the image of $\otimes F_i$ via $\alpha_r$. Since it is symmetric in the 
sets $\ux,\uy$, it can be thought of as an element of $S_2(S_{re})$. 
\end{itemize} 

Suppose now that $\dim V = 2$. We will show by induction on $r$ that 
$\pi_p \circ \alpha_r$ is not identically zero. 

\subsection{Case $r=2$.} 

We now specialize the $F_i$. 
Let $F_i = l_i^d$, where $l_i(x_0,x_1)$ are linear forms. 
Then $\alpha_r(\otimes F_i) = Q(\ux)^eQ(\uy)^e$,
where $Q = \prod l_i$. 

Introduce the Omega operator 
\[ \Omega = 
\frac{\partial^2}{\partial x_0 \, \partial y_1} 
-\frac{\partial^2}{\partial y_0 \, \partial x_1}. \]
The projection $\pi_p$ corresponds to applying 
$\Omega^{2p}$ and substituting $\uy:=\ux$, which is the same as 
taking the $2p$-th transvectant of $Q^e$ with itself. 
Hence we are reduced to showing the following statement: 

\smallskip

\begin{Lemma} \sl 
If $Q$ is the generic binary quadratic, then 
\[ (Q^e,Q^e)_{2p} \neq 0, \quad \text{for $0 \le p \le e$}. \] 
\label{lemma.A} \end{Lemma}

\demo See Section \ref{Feynman}. \qed 

\smallskip 

\subsection{The induction step} 
For the transition from $r$ to $r+1$, consider the commutative diagram 
\[ \diagram 
S_r(S_d) \otimes S_d \dto \rto^{\alpha_r \otimes 1} & 
H^0(\O_X(r)) \otimes S_d \dto^{u_r} \\ 
S_{r+1}(S_d) \rto_{\alpha_{r+1}} & H^0(\O_X(r+1)) \enddiagram \] 
Assume that $\alpha_r$ (and hence $\alpha_r \otimes 1$) is surjective. If we show that $u_r$ is 
surjective, then it will follow that $\alpha_{r+1}$ is 
surjective. We need to understand the action of $u_r$ on the 
summands of the decomposition (\ref{oxr}). The map 
\[ u_r^{(p,p')}: S_{rd-4p} \otimes S_d \lra S_{(r+1)d-4p'} \] 
is defined as the composite 
\[ \begin{aligned} 
{} & S_{rd-4p} \otimes S_d \ra (S_{re} \otimes S_{re}) \otimes S_d 
\ra (S_{re} \otimes S_{re}) \otimes (S_e \otimes S_e) \ra \\
& S_{(r+1)e} \otimes S_{(r+1)e} \ra S_{(r+1)d-4p'}. 
\end{aligned} \] 
Let $C \in S_{rd-4p},D \in S_d$. Following the component maps, 
we will get a recipe for calculating the image 
of $C \otimes D$ via $u_r^{(p,p')}$. Let 
\[ \Gamma_a = \sum \limits_{i=0}^{re} \binom{re}{i} \, a_i \, x_0^{re-i} x_1^i, \quad
\Gamma_b = \sum \limits_{i=0}^{re} \binom{re}{i} \, b_i \, x_0^{re-i} x_1^i, \]
be two {\sl generic} forms of degree $re$. (That is to say, the 
$a_i,b_i$ are thought of as independent indeterminates.) 
\begin{itemize} 
\item 
Let $T_1 = (\Gamma_a,\Gamma_b)_{2p}$ and 
$T_2 = (C,T_1)_{rd-4p}$. Then $T_2$ does not involve $x_i$. 
\item 
Obtain $T_3$ by making the substitutions 
\[ a_i = x_1^{re-i} (-x_0)^i, \quad 
   b_i = y_1^{re-i} (-y_0)^i, \] 
in $T_2$. 
\item 
Let 
\[ T_4 = (y_0 \frac{\partial}{\partial x_0} + 
          y_1 \frac{\partial}{\partial x_1})^{e} \, D, \] 
and $T_5 = T_3 \, T_4$. 
\item 
Let $T_6 = \Omega^{2p'} \, T_5$. 
Finally $u_r^{(p,p')}(C \otimes D)$ is obtained by substituting 
$x_0,x_1$ for $y_0,y_1$ in $T_6$. 
\end{itemize} 
Hence it is enough to show the following: 

For $p'$ in the range $0 \le p' \le \frac{(r+1)e}{2}$, there exists a $p$ such 
that $u_r^{(p,p')}(C \otimes D)$ is nonzero for some forms $C,D$ of 
degrees $rd-4p,d$ respectively. 

We translate this statement into the symbolic calculus of classical invariant 
theory (see \cite{GrYo}). Introduce symbolic letters $c,d$, and 
let $c_{\ux}$ stand for $c_0 x_0 + c_1x_1$ etc. 
Write $\omega = x_0y_1 - x_1y_0$. Then the claim becomes

\smallskip 

\begin{Lemma} \sl 
Given $r \ge 2$ and $0 \le p' \le \frac{(r+1)e}{2}$, there exists 
a $p$ in the range $0 \le p \le \frac{re}{2}$, such that the expression 
\[ 
\{ 
\Omega^{2p'}( \omega^{2p} \, c_{\ux}^{re-2p} \, c_{\uy}^{re-2p} d_{\ux}^{\,e} \, 
d_{\uy}^{\,e}) \}|_{\uy:=\ux} \] 
is nonzero. 
\label{lemma.B} \end{Lemma}

\demo See Section \ref{Feynman}. \qed 

\medskip 

At this point, modulo Lemmata \ref{lemma.A} and \ref{lemma.B}, 
the proof of the main theorem is complete. 

\begin{Remark} \rm 
It would be unnecessary to make an inductive argument if we could prove the 
following statement: 

{\sl Assume $r \ge 2$, and let $R$ denote the generic binary form of 
degree $r$. Then $(R^e,R^e)_{2p} \neq 0$ 
for $0 \le p \le \lfloor \frac{re}{2} \rfloor$.}

However we do not see how to do this. 
\end{Remark} 

\section{The combinatorics of Feynman diagrams} \label{Feynman}

In this section we will complete the proof of the key
Proposition \ref{prop.nonzero}, by proving Lemmata \ref{lemma.A} and \ref{lemma.B}.
Although there might be a simple geometric or representation
theoretic argument allowing the derivation of these lemmas,
we were unable to find one, and relied instead
on explicit combinatorial computation.
This is an instance of what one might call the combinatorics
of invariants of binary forms which
were at the heart of {\em classical}
invariant theory.
Although neglected for the last century, it is a fascinating subject with
ramifications in many fields of current mathematical and physical interest
like the theory of angular momentum~\cite{Biedenharn1, Biedenharn2},
classical hypergeometric series~\cite{Gustafson},
the spin network approach to quantum gravity~\cite{Penrose, Rovelli},
as well as knot and 3-manifold invariants~\cite{Carter}.

As for lemma \ref{lemma.A}, we will prove that
\be
\left.
\Omega_{xy}^{2p}
Q(x)^e Q(y)^e \right|_{y=x}
=
\cN_{e,p}^{\rm I}
Q(x)^{2e-2p}
(-\De)^{p}
\ee
where $\De$ is the discriminant of the quadratic form $Q$ and
$\cN_{e,p}^{\rm I}$ is a strictly positive numerical constant.
We will give two proofs of this result.
The first is a combinatorially explicit calculation with Feynman diagrams
which explains {\em why} $\cN_{e,p}^{\rm I}>0$.
The second is perhaps less transparent, but it
allows the exact computation of the constant $\cN_{e,p}^{\rm I}$. 
Apart from a harmless normalisation factor, it is 
a special value of Wigner's $3j$-symbol (or Clebsch-Gordan coefficient,
see~\cite{Gustafson}) which can be computed using
Dixon's summation theorem for the ${}_3 F_2$ hypergeometric series.
The comparison of both methods yields an interesting
formula for the weighted enumeration of a class of bipartite graphs
which have vertex degree at most two.
Note that  Proposition \ref{prop.nonzero} can also be considered as a statement
concerning a sum over bipartite graphs with vertex degree bounded
by r, the degree of the binary form $Q$.
However a direct combinatorial approach seems very difficult at this point.
The representation theoretic arguments involved in our inductive
proof of Proposition \ref{prop.nonzero} and its reduction to 
Lemmata \ref{lemma.A} and \ref{lemma.B}, can be
credited for taming a significant part of the combinatorial
complexity of such sums over graphs.

The proof of Lemma \ref{lemma.B}, uses Feynman diagrammatic generating
function techniques which are implicit in the work
of J.Schwinger~\cite{Schwinger} and its reformulation by
V. Bargmann~\cite{Bargmann}.
This allows us to prove that
\be
\left.
\Omega_{xy}^{2p'} \, 
c_x^{re-2p} \, c_y^{re-2p} \, d_x^e \, d_y^e
\right|_{y=x}=
\cN_{r,e,p',p}^{\rm II}\times \, 
(cd)^{p'-p} c_x^{2(re-p'-p)} \, d_x^{2(e-p'+p)},
\ee
where
$\cN_{r,e,p',p}^{\rm II}$ is a numerical constant that we compute
explicitly.
An easy and tempting shortcut at this
point would have been to use analysis, akin to what Bargmann did in~\cite{Bargmann}.
This would however obscure the fact that what is at play
is purely combinatorial algebra, with no real need for transcendental
methods.

\subsection{First proof of Lemma \ref{lemma.A}}

The following presentation is semi-formal yet completely rigorous.
The reader who needs a stricly formal exposition of Feynman diagrams
and their rigorous mathematical use should consult~\cite{Abdesselam}
(see also~\cite{Fiorenza}).
For the present purposes, let us simply say that a Feynman diagram
is essentially the combinatorial data needed to encode
a complex tensorial expression built from a predefined collection
of elementary tensors, exclusively using contraction of
tensor indices. The word ``tensor'' here is used
as meaning the multidimensional analog of a matrix which
is therefore basis dependent.
Coordinates are needed in order to {\em state}
the necessary definitions, but are almost never actually {\em used}
in the computations.
Here the basic tensors are
\be
x=\lp
\begin{array}{c}
x_1\\
x_2
\end{array}
\rp
\ \ \ ,\ \ \ 
y=\lp
\begin{array}{c}
y_1\\
y_2
\end{array}
\rp
\ee
made of formal indeterminates,
and the two matrices $Q$ and $\ep$
in $M_2(\CC)$.
$Q$ is symmetric and gives the quadratic form
$Q(x)=x^{\rm T}Qx$,
\be
\ep\eqdef\lp
\begin{array}{cc}
0 & 1 \\
-1 & 0
\end{array}
\rp 
\ee
is antisymmetric and defines the symbolic brackets
as well as Cayley's Omega operator.
We also need the vectors of diffential
operators
\be
\partial_x=\lp
\begin{array}{c}
\frac{\partial}{\partial x_1}\\
\frac{\partial}{\partial x_2}
\end{array}
\rp
\ \ \ ,\ \ \ 
\partial_y=\lp
\begin{array}{c}
\frac{\partial}{\partial y_1}\\
\frac{\partial}{\partial y_2}
\end{array}
\rp
\ee
We now introduce a graphical notation for the
entries of these elementary tensors (indices belong to the set $\{1,2\}$),

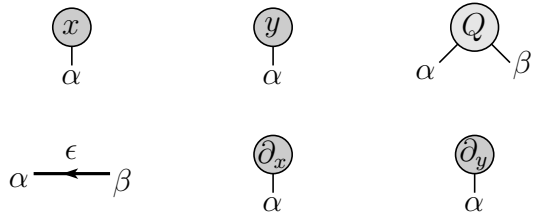
\[
\figplace{dessin1}{0 in}{0 in}
\figplace{dessin2}{0 in}{0 in}
\figplace{dessin3}{0 in}{0 in}
\]
\[
\figplace{dessin4}{0 in}{0 in}
\figplace{dessin5}{0 in}{0 in}
\figplace{dessin6}{0 in}{0 in}
\]

Now to any diagram obtained by assembling
any number of these elementary pieces by gluing
pairs of index-bearing lines,
one associates an expression, called the {\em amplitude}
of the diagram.
For example

\be
\figplace{dessin7}{0 in}{-0.28 in}
\eqdef
\sum_{\al,\beta=1}^2
x_\al Q_{\al\beta} x_\beta
=Q(x)
\ee
the quadratic form itself.
also
\bea
\figplace{dessin8}{0 in}{-0.2 in} & \eqdef & \sum_{\al,\beta,\ga,\de=1}^2
Q_{\al\beta}\ep_{\al\ga}\ep_{\beta\de} Q_{\ga\de} \\
 & = & 2(Q_{11}Q_{22}-Q_{12}Q_{21}) \\
 & = & 2 \det(Q)
\eea
Whenever we write a diagram inside an equation what is meant
is the amplitude of the diagram.
Now we will use the fact $Q_{\al\beta}$ has an {\em inner
structure}. This is related to the notion of combinatorial plethysm
(see~\cite{Joyal}).
Indeed, since $\CC$ is algebraically closed, one can
factor $Q$ as
$Q(x)=R_1(x) R_2(x)$
where 
\be
R_1=\lp
\begin{array}{c}
R_{1,1}\\
R_{1,2}
\end{array}
\rp
\ \ \ ,\ \ \ 
R_2=\lp
\begin{array}{c}
R_{2,1}\\
R_{2,2}
\end{array}
\rp
\ee
in $\CC^2$, are dual to the homogenous roots of $Q$.
Now for any indices $\al$ and $\beta$
\bea
Q_{\al\beta} & = & \frac{1}{2}
\frac{\partial^2}{\partial x_\al
\partial x_\beta}
Q(x) \\
 & = & \frac{1}{2}
\lp
R_{1,\al} R_{2,\beta} +
R_{1,\beta} R_{2,\al}
\rp
\eea
which we write more suggestively as

\be
\figplace{dessin3}{0 in}{-0.25 in}
=
\figplace{dessin9}{0 in}{-0.25 in}
\ee
\be
=\frac{1}{2}
\figplace{dessin10}{0 in}{-0.5 in}
+\frac{1}{2}
\figplace{dessin11}{0 in}{-0.5 in}
\label{Decomp}
\ee

This implies, for instance that

\be
\figplace{dessin8}{0 in}{-0.2 in}
=
\frac{1}{4}
{{
\figplace{dessin12}{0 in}{-0.1 in}
}\atop{
\figplace{dessin13}{0 in}{-0.1 in}
}}
+
\frac{1}{4}
{{
\figplace{dessin13}{0 in}{-0.1 in}
}\atop{
\figplace{dessin12}{0 in}{-0.1 in}
}}
\ee
since reversing the direction of an $\ep$ arrow
produces a minus sign, and therefore
\be
\figplace{dessin12bis}{0 in}{-0.1 in}
=
\figplace{dessin13bis}{0 in}{-0.1 in}
=0
\ee 
As a result
\be
\figplace{dessin8}{0 in}{-0.2 in}
= -\frac{1}{2}\De
\ee
where
\be
\De\eqdef
\lp\figplace{dessin12}{0 in}{-0.1 in}\rp^2
\ee
is the discriminant of $Q$.

Now the quantity we are interested in is
\be
\left.
\Omega_{xy}^{2p}
Q(x)^e Q(y)^e \right|_{y=x}
=F(x,x)
\ee
where
\be
F(x,y)
=
\lp\figplace{dessin14}{0 in}{-0.09 in}
\rp^{2p}
\lp\figplace{dessin7}{0 in}{-0.26 in}
\rp^e
\lp\figplace{dessin15}{0 in}{-0.26 in}
\rp^e
\label{Fofxy}
\ee
Summing over where exactly the derivatives act, via Leibnitz's rule,
generates a sum over
Feynman diagrams which, once we let $y=x$, condenses
into the following sum over vertex-labelled
bipartite multigraphs.
\be
F(x,y)
=
\sum_G w_G \cA_G
\label{SumG}
\ee
Here, the graph $G$ can be seen as a matrix $(m_{ij})$
in $\NN^{L\times R}$ where $L$ and $R$
are fixed sets of cardinality $e$ labelling
the $Q(x)$ and $Q(y)$ factors in (\ref{Fofxy})
respectively.
$w_G$ is a combinatorial weight and
$\cA_G$ is the amplitude of a Feynman diagram
encoded by the graph $G$.
The latter has to satisfy the following conditions
\[
\sum_{i\in L, j\in R} m_{ij}=2p
\]
\[
\forall i\in L,\ 
l_i\eqdef\sum_{j\in R} m_{ij}\le 2
\]
and
\be
\forall j\in R,\ 
c_j\eqdef\sum_{i\in L} m_{ij}\le 2
\label{ConditionG}
\ee
The combinatorial weight is easily seen to be
\be
w_G=
\frac{(2p)! 2^{2e}}
{\prod_{i,j} (m_{ij})! \times 
\prod_i (2-l_i)! \times
\prod_j (2-c_j)!}
\ee

The amplitude $\cA_G$ factors over the connected components
of $G$. These components are of four possible types:
cycles containing an even number of $\ep$ arrows
of alternating direction,
chains with both endpoints in $L$,
chains with both endpoints in $R$ and finally
chains with one endpoint in $L$ and another in $R$.
However, the last type of connected component
gives a zero contribution.
Indeed, such a chain contains an odd number of $\ep$
arrows and therefore its amplitude changes sign
if one reverses all the orientations of the arrows,
an operation which, followed by a 180 degrees rotation, puts
the chain back in its original form; because we have already set $y=x$.
For exemple

\be
\figplace{dessin16}{0 in}{-0.45 in}
=
\figplace{dessin17}{0 in}{-0.18 in}
\ee
\be
=\ -\figplace{dessin18}{0 in}{-0.18 in}
\ee
\be
=\ -\figplace{dessin17}{0 in}{-0.18 in}
\ee
\be
=\ 0
\ee

Now to calculate the amplitudes of the other three kinds of components,
one needs to use the inner structure of $Q$.
Namely for each cycle of even length $2m$, incorporating
the decomposition (\ref{Decomp})
at each $Q$ vertex, produces a sum of $2^{2m}$ terms
all of which vanish except for two of them.
Indeed, once one choses the precise connections between
the ``inner'' and ``outer'' part of what was a particular $Q$ vertex,
the connections for the remaining vertices are forced, if
one wants to avoid the appearance of the vanishing factors

\[
\figplace{dessin12bis}{0 in}{-0.09 in}
\ \ \ {\rm and}\ \ \   
\figplace{dessin13bis}{0 in}{-0.09 in}
\]

Besides, the alternating pattern for the orientations of the $\ep$ arrows
makes it so that we collect an equal number $m$ of
\[
\figplace{dessin12}{0 in}{-0.09 in}
\ \ \ {\rm and}\ \ \   
\figplace{dessin13}{0 in}{-0.09 in}
\]
factors.
As a result, the amplitude of the cycle is exactly
$2^{1-2m} (-\De)^m$.

Likewise, a chain with both endpoints in $L$, or both endpoints in $R$,
and with a necessarily even number $2m$ of $\ep$
arrows (and thus $2m+1$ $Q$ vertices) gives as an amplitude

\[
\frac{2}{2^{2m+1}}
\figplace{dessin19}{0 in}{-0.09 in}
\lp
\figplace{dessin13}{0 in}{-0.09 in}
\figplace{dessin12}{0 in}{-0.09 in}
\rp^m
\figplace{dessin20}{0 in}{-0.09 in}
\]
\be
=
2^{-2m} (-\De)^m Q(x)
\ee

Therefore
an easy count shows that the amplitude of a bipartite graph $G$ in
(\ref{SumG})
is
\be
\cA_G=
2^{\cC(G)-2p}
Q(x)^{2e-2p} (-\De)^p
\ee
where $\cC(G)$ is the number of cycles in $G$.
Finally,
\be
F(x,x)=\cN_{e,p}^{\rm I}
Q(x)^{2e-2p} (-\De)^p
\ee
where
\be
\cN_{e,p}^{\rm I}\eqdef
\sum_{G}
\frac{(2p)! 2^{2e-2p+\cC(G)}}
{\prod_{i,j} (m_{ij})! \times 
\prod_i (2-l_i)! \times
\prod_j (2-c_j)!}
\label{NIdef}
\ee
and the last sum is over all graphs
$G=(m_{ij})$ satisfying the three constraints (\ref{ConditionG})
{\em and} the additional condition
of not having any connected component which is a chain starting
in $R$ and ending in $L$.
It is easy to see that given $e\ge 1$ and $p$, $0\le p\le e$,
there always exists such graphs $G$, i.e.
$\cN_{e,p}^{\rm I}>0$ which proves Lemma \ref{lemma.A}.

\subsection{Second proof of Lemma \ref{lemma.A}}

We specialize the quadratic form to $Q(x)=x_1 x_2$, for which
$\De=1$.
We use
\be
\Om_{xy}^{2p}
=
\sum_{i=0}^{2p}
(-1)^i
\lp
\begin{array}{c}
2p \\
i
\end{array}
\rp
\frac{\partial^{2p}}{\partial x_1^{2p-i} \partial x_2^{i}}
\frac{\partial^{2p}}{\partial y_1^{i} \partial y_2^{2p-i}}
\ee
to obtain
\[
\left.
\Om_{xy}^{2p}
Q(x)Q(y)
\right|_{y=x}
=
\sum_{i=\max(0,2p-e)}^{\min(2p,e)}
(-1)^i
\lp
\begin{array}{c}
2p \\
i
\end{array}
\rp\times
\]
\be
\frac{e!^4}{(e-2p+i)!^2(e-i)!^2}
\times x_1^{2e-2p} x_2^{2e-2p}
\ee
\be
=\cN_{e,p}^{\rm I} Q(x)^{2e-2p} (-\De)^p
\label{NIred}
\ee
with
\be
\cN_{e,p}^{\rm I}
=
\sum_{i=\max(0,2p-e)}^{\min(2p,e)}
(-1)^{p+i}
\lp
\begin{array}{c}
2p \\
i
\end{array}
\rp
\frac{e!^4}{(e-2p+i)!^2(e-i)!^2}
\label{NIdef2}
\ee
By $GL_2(\CC)$ change of coordinate and density of quadratic forms
with nonzero discriminants, (\ref{NIred}) is valid for any
quadratic form $Q$.

Let us suppose $e\ge 2p\ge 0$. It is then easy to rewrite
$\cN_{e,p}^{\rm I}$, using Pochammer's symbol
$(a)_n\eqdef a(a+1)\cdots (a+n-1)$,
as
\bea
\cN_{e,p}^{\rm I} & = & \frac{(-1)^p e!^2}{(e-2p)!^2}
\sum_{i\ge 0}
\frac{(-2p)i (-e)_i (-e)_i}
{i! (e-2p+1)_i (e-2p+1)_i} \\
 & = & \frac{(-1)^p e!^2}{(e-2p)!^2}
{ }_3 F_2
\left[
\begin{array}{c}
-2p,-e,-e \\
e-2p+1, e-2p+1
\end{array}
; 1
\right]
\eea
The terminating classical hypergeometric series that
appears
in the last formula, is of the form
\[
{ }_3 F_2
\left[
\begin{array}{c}
a,b,c \\
1+a-b, 1+a-c
\end{array}
; 1
\right]
\]
and can therefore be evaluated thanks to Dixon's summation
theorem (see~\cite{Slater}) :
\bea
\lefteqn{
{ }_3 F_2
\left[
\begin{array}{c}
a,b,c \\
1+a-b, 1+a-c
\end{array}
; 1
\right]=} & & \nonumber\\
 & &
\frac{\Ga(1+\frac{1}{2}a)\Ga(1+\frac{1}{2}a-b-c)
\Ga(1+a-b)\Ga(1+a-c)}
{\Ga(1+a)\Ga(1+a-b-c)
\Ga(1+\frac{1}{2}a-b)\Ga(1+\frac{1}{2}a-c)}
\eea
which is valid in the domain of analyticity $\Re({1+\frac{1}{2}a-b-c})>0$.
Here we want to take $a=-2p$ and $b=c=-e$; one therefore has to be
careful with the
$\frac{\Ga(1+\frac{1}{2}a)}{\Ga(1+a)}$ factor and rewrite it
as
\be
\frac{\pi}{\Ga(-\frac{a}{2})\sin(-\frac{\pi a}{2})}
\times\frac{\Ga(-a)\sin(-\pi a)}{\pi}=
\cos(\frac{\pi a}{2})\frac{\Ga(-a+1)}{\Ga(-\frac{a}{2}+1)}
\ee
The end result is
\be
\cN_{e,p}^{\rm I}
=
\frac{(2p)! (2e-p)! e!^2}
{p!(2e-2p)!(e-p)!^2}
\ee
if $0\le p\le \frac{e}{2}$.

Now if $\frac{e}{2}\le p\le e$,
then by setting $i=2p-e+j$
in (\ref{NIdef2}) and again writing the resulting sum over 
$j$ as a terminating hypergeometric series
one gets
\be
\cN_{e,p}^{\rm I}
=
\frac{(-1)^{p+e}(2p)! e!^3}
{(2p-e)!(2e-2p)!^2}
{}_3 F_2
\left[
\begin{array}{c}
-2e+2p, -2e+2p,-e \\
1, 2p-e+1
\end{array}
; 1
\right]
\ee
The same method using Dixon's theorem gives
\be
\cN_{e,p}^{\rm I}
=
\frac{(2p)!(3e-3p)! e!^3}
{(2p-e)!(2e-2p)!^2(e-p)!^3}
\ee
if $\frac{e}{2}\le p\le e$.

This again shows that, in either case, $\cN_{e,p}^{\rm I}>0$
and completes our second proof of Lemma \ref{lemma.A}.
This also gives a closed form evaluation of the sum
in (\ref{NIdef})

\subsection{Proof of Lemma \ref{lemma.B}}

Let $r$, $e$, $p'$ and $p$ be integers satisfying
$r\ge 2$, $e\ge 1$, $0\le 2p'\le (r+1)e$ and $0\le 2p\le re$.
Let
\be
c=\lp
\begin{array}{c}
c_1\\
c_2
\end{array}
\rp
\ \ \ ,\ \ \ 
d=\lp
\begin{array}{c}
d_1\\
d_2
\end{array}
\rp
\ee
be two elements of $\CC^2$
and
\be
x=\lp
\begin{array}{c}
x_1\\
x_2
\end{array}
\rp
\ \ \ ,\ \ \ 
y=\lp
\begin{array}{c}
y_1\\
y_2
\end{array}
\rp
\ee
be two vectors of indeterminates.
The quantity we would like to compute is
\be
G(x)\eqdef
\left.
\Om_{xy}^{2p'}
(xy)^{2p}
c_x^{re-2p}c_y^{re-2p} d_x^e d_y^e
\right|_{y=x}
\ee
or, in matrix notation,
\be
G(x)=
\left.
[\partial_x^{\rm T}\ep\partial_y]^{2p'}
(x^{\rm T}\ep y)^{2p}
(x^{\rm T}c c^{\rm T} y)^{re-2p}
(x^{\rm T}d d^{\rm T} y)^{e}
\right|_{y=x}
\ee
We now introduce two new vectors of inderminates
\be
\ph=\lp
\begin{array}{c}
\ph_1\\
\ph_2
\end{array}
\rp
\ \ \ ,\ \ \ 
{\Br\ph}=\lp
\begin{array}{c}
{\Br\ph}_1\\
{\Br\ph}_2
\end{array}
\rp
\ee
It is then easy to see that
\be
G(x)=(2p')!(2p)!(re-2p)!e!
\frac{[\partial_\ph^{\rm T}\ep\partial_{\Br\ph}]^{2p'}}
{(2p')!}
\frac{[(\ph+x)^{\rm T}\ep({\Br\ph}+x)]^{2p}}
{(2p)!}
\nonumber
\ee
\be
\left.
\times
\frac{[({\Br\ph}+x)^{\rm T}c c^{\rm T}(\ph+x)]^{re-2p}}
{(re-2p)!}
\times
\frac{[({\Br\ph}+x)^{\rm T}d d^{\rm T}(\ph+x)]^{e}}
{e!}
\right|_{{\ph=0}\atop{{\Br\ph}=0}}
\ee
i.e.
\be
G(x)=(2p')!(2p)!(re-2p)!e!
[h^{2p'}u^{2p}v^{re-2p}w^e]\cZ
\ee
where $[h^{2p'}u^{2p}v^{re-2p}w^e]\cZ\in\CC[[x_1,x_2]]$
denotes the coefficient of the monomial
$h^{2p'}u^{2p}v^{re-2p}w^e$ in
$\cZ\in\CC[[x_1,x_2,h,u,v,w]]$, the generating function
defined by
\be
\cZ\eqdef
\left.
\left\{
\sum_{n\ge 0} \frac{h^n}{n!}
[\partial_\ph^{\rm T}\ep\partial_{\Br\ph}]^{n}
e^S
\right\}
\right|_{{\ph=0}\atop{{\Br\ph}=0}}
\ee
where $S\in\CC[[\ph_1\ph_2,{\Br\ph}_1,{\Br\ph}_2,x_1,x_2,h,u,v,w]]$
is given by
\be
S\eqdef
({\Br\ph}+x)^{\rm T}(-u\ep+M)(\ph+x)
\ee
and
\be
M\eqdef
v c c^{\rm T}+ w d d^{\rm T}\in
M_2(\CC[[v,w]]).
\ee
Note that there is no problem of convergence since we work over
rings of formal power series with their usual topology.
With obvious notations, one can rewrite
$\cZ$ as
\be
\cZ=
\left.
\exp(h\partial_\ph^{\rm T}\ep\partial_{\Br\ph})
\exp({\Br\ph}^{\rm T} A\ph+J^{\rm T}\ph+{\Br\ph}^{\rm T}K
+S_0)
\right|_{{\ph=0}\atop{{\Br\ph}=0}}
\ee
with $A\eqdef -u\ep+M$, $J^{\rm T}\eqdef
-ux^{\rm T}\ep+x^{\rm T}M$,
$K=-u\ep x+Mx$ and
$S_0\eqdef v (x^{\rm T}c)^2+w (x^{\rm T}d)^2$.
Therefore $\cZ=e^{S_0}{\til\cZ}$ with
\be
{\til\cZ}\eqdef
\left.
\exp(h\partial_\ph^{\rm T}\ep\partial_{\Br\ph})
\exp({\Br\ph}^{\rm T} A\ph+J^{\rm T}\ph+{\Br\ph}^{\rm T}K)
\right|_{{\ph=0}\atop{{\Br\ph}=0}}
\ee

${\til\cZ}$ can now be expressed as a sum over
Feynman diagrams built, like in Section 5.1, from the following pieces
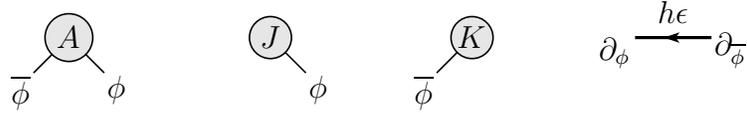
\[
\figplace{dessin21}{0 in}{0 in}
\figplace{dessin22}{0 in}{0 in}
\figplace{dessin23}{0 in}{0 in}
\figplace{dessin24}{0 in}{0 in}
\]
by plugging the $\partial_\ph$'s
onto the $\ph$'s and 
the $\partial_{{\Br\ph}}$ onto the ${\Br\ph}$'s, in all possible ways.
More precisely, given any finite set $E$, we define a Feynman diagram
on $E$ as any sextuple $\cF=(E_\ph,E_{\Br\ph},
\pi_A,\pi_J,\pi_K,\cC)$
where $E_\ph$, $E_{\Br\ph}$ are subsets of $E$,
and $\pi_A$, $\pi_J$, $\pi_K$ are, each, sets of subsets of $E$
and $\cC$ is a map $E_{\Br\ph}\rightarrow E_\ph$, satisfying the following
axioms
\begin{itemize}
\item
$E_\ph$ and $E_{\Br\ph}$
have equal cardinality and they form a two set partition of $E$.
\item
The union of the elements in $\pi_A$,
that of the elements of $\pi_J$, and likewise for $\pi_K$
form a three set partition of $E$.
\item
$\cC$ is bijective.
\item
Every element of $\pi_A$ has two elements,
one in $E_\ph$ and one in $E_{\Br\ph}$. 
\item
Every element of $\pi_J$ has only one element
which is in $E_\ph$.
\item
Every element of $\pi_K$ has only one element
which is in $E_{\Br\ph}$.
\end{itemize}

The set of Feynman diagrams on $E$
is denoted by ${\mathsf{Fey}}(E)$.
The point is that the association
$E\rightarrow {\mathsf{Fey}}(E)$
defines an endofunctor of the groupoid category of finite
sets with bijections~\cite{Joyal, Abdesselam, Fiorenza},
since, given a Feynman diagram $\cF$ on $E$ and a bijective
map $\si:E\rightarrow E'$, there is a natural
way to transport $\cF$ along $\si$ to
obtain a Feynman diagram $\cF'={\mathsf{Fey}}(\si)(\cF)$
on $E'$. Let us take for example
$E=\{1,2,\ldots,8\}$, $E_\ph=\{1,2,3,4\}$,
$E_{\Br\ph}=\{5,6,7,8\}$,
$\pi_A=\{ \{2,6\},\{3,7\}, \{4,8,\} \}$,
$\pi_J=\{\{1\}\}$, $\pi_K=\{\{5\}\}$,
and $\cC$ given by $\cC(5)=1$, $\cC(6)=3$,
$\cC(7)=4$ and $\cC(8)=2$.
This corresponds to the diagram
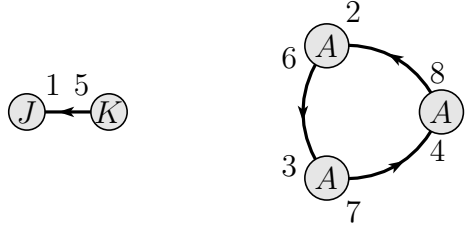
\[
\figplace{dessin25}{0 in}{0.8 in}
\figplace{dessin26}{0 in}{0 in}
\]
where we put the elements of $E$ next to the corresponding half-line.
The amplitude of such a pair $(E,\cF)$ is in this example
\be
\cA(E,\cF)=(J^{\rm T}(h\ep)K)\times tr([h\ep A]^3)
\ee
Note that there is a natural equivalence relation between pairs
of finite sets equiped with a Feynman diagram. It is given by
letting
$(E,\cF)\sim(E',\cF')$ if and only if there exists a bijection
$\si:E\rightarrow E'$ such that $\cF'={\mathsf{Fey}}(\si)(\cF)$.
One also has the notion of automorphism group $Aut(E,\cF)$
of a pair $(E,\cF)$ which is the set of bijections $\si:E\rightarrow E$
such that ${\mathsf{Fey}}(\si)(\cF)=\cF$.
Now one can check that
\be
{\til\cZ}=\sum_{[E,\cF]}\frac{\cA(E,\cF)}{|Aut(E,\cF)|}
\ee
where the sum is over equivalence classes of pairs
$(E,\cF)$, $\cA(E,\cF)$
is the amplitude, and $|Aut(E,\cF)|$ the cardinality of
the automorphism group of any representative in the class $[E,\cF]$.
We leave it to the reader to check (otherwise see~\cite{Abdesselam,
Fiorenza}) that
\be
\log {\til\cZ}=
\sum_{[E,\cF]\ \rm connected}
\frac{\cA(E,\cF)}{|Aut(E,\cF)|}
\ee
\be
= \sum_{n\ge 1} \frac{1}{n} tr((h\ep A)^n)
+
\sum_{n\ge 0}
J^{\rm T}(h\ep A)^n(h\ep)K
\ee
since the only connected diagrams are pure $A$ cycles
or $A$ chains joining a $J$ to a $K$ vertex.
As a result
\be
{\til\cZ}=\frac{1}{\det(I-h\ep A)}
\exp(J^{\rm T}(I-h\ep A)^{-1}(h\ep)K)
\ee
After straightforward but tedious computations with 2 by 2
matrices, which we spare the reader,
this implies
\be
\cZ=\frac{1}{(1-hu)^2+h^2vw(c^{\rm T}\ep d)^2}
\exp\lp
\frac{v(x^{\rm T}c)^2+w(x^{\rm T}d)^2}{(1-hu)^2+h^2vw(c^{\rm T}\ep d)^2}
\rp
\ee
or in classical notation
\be
\cZ=\frac{1}{(1-hu)^2+h^2vw(cd)^2}
\exp\lp
\frac{v c_x^2+w d_x^2}{(1-hu)^2+h^2vw(cd)^2}
\rp
\ee
We now expand
\be
\cZ=\sum_{\mu\ge 0}\frac{1}{\mu!}
(v c_x^2+w d_x^2)^{\mu}
\lp
1-2hu+h^2u^2+h^2vw(cd)^2
\rp^{-(\mu+1)}
\ee
\be
=\sum_{\mu,\nu\ge 0}
\frac{(-1)^\nu (\mu+\nu)!}{\mu!^2\nu!}
(v c_x^2+w d_x^2)^{\mu}
\lp
-2hu+h^2u^2+h^2vw(cd)^2
\rp^{\nu}
\ee
\be
=\sum_{{{m,n}\atop{\al,\beta,\ga}}\ge 0}
\frac{(-1)^{\al+\beta+\ga}(m+n+\al+\beta+\ga)!}
{(m+n)!m!n!\al!\beta!\ga!}
\times
\ee
\be
(v c_x^2)^m (w d_x^2)^n (-2hu)^\al
(h^2u^2)^\beta {\lp h^2vw(cd)^2\rp}^\ga
\ee
\be
=\sum_{{{m,n}\atop{\al,\beta,\ga}}\ge 0}
\frac{(-1)^{\beta+\ga}(m+n+\al+\beta+\ga)! 2^\al}
{(m+n)!m!n!\al!\beta!\ga!}
h^{\al+2\beta+\ga}
u^{\al+2\beta}v^{m+\ga}w^{n+\ga}
\times
\ee
\be
c_x^{2m} d_x^{2n} (cd)^{2\ga}
\ee
The coefficient of $h^{2p'}u^{2p} v^{re-2p} w^{e}$
is a sum over the single index $\beta$, $0\le \beta\le p$,
as a result of solving for $\al=2p-2\beta$, $\ga=p'-p$,
$m=re-p'-p$,
$n=e-p'+p$.
Therefore
\be
G(x)=\cN_{r,e,p',p}^{\rm II}
c_x^{2(re-p'-p)} d_x^{2(e-p'+p)} (cd)^{2(p'-p)}
\ee
where
\be
\cN_{r,e,p',p}^{\rm II}=
\bbbone_{\left\{
{p'-p\ge 0,e-p'+p\ge 0}\atop{re-p'-p\ge 0}
\right\}}
\times
\ee
\be
\frac{(-1)^{p'-p} (2p)!(2p')!(re-2p)!e!}
{(p'-p)!(e-p'+p)!(re-p'-p)!((r+1)e-2p')!}
\times \cJ_{s,p}
\ee
where
$\bbbone_{\{\cdots\}}$ denotes the characteristic function
of the condition between braces and
\be
\cJ_{s,p}\eqdef
\sum_{\beta=0}^p \frac{(-1)^\beta 2^{2p-2\beta} (s+2p-\beta)!}
{(2p-2\beta)!\beta!}
\ee
with $s\eqdef
(r+1)e-p'-p\ge e$
when the characteristic function is nonzero.
Note that $\cJ_{s,p}$ can be rewritten as a Gauss hypergeometric series
and can be summed by the Chu-Vandermonde theorem (see~\cite{Slater} for
instance) 
\be
\cJ_{s,p}=\frac{(s+p)!(s+\frac{3}{2})_p}{p!(\frac{1}{2})_p}
\ee
As a result, the characteristic function alone dictates whether 
$\cN_{r,e,p',p}^{\rm II}$ vanishes or not. 
Now for $r\ge 2$, $e\ge 1$ and
$0\le p'\le \frac{(r+1)e}{2}$ it is easy to see that one can allways
find an integer $p$
with $0\le p\le \frac{re}{2}$,
$p'-p\ge 0$, $e-p'+p\ge$ and $re-p'-p\ge 0$.
Indeed take $p=p'$ if $0\le p'\le \frac{re}{2}$,
and otherwise take $p=p'-e$ if $\frac{re}{2}\le p'\le \frac{(r+1)e}{2}$.
This completes the proof of Lemma \ref{lemma.B}.

\section{A note on terminology and history}

The simultaneous relevance to our approach of the literature from many fields of
mathematics and physics requires the following {\it mise au point}.
Firstly, our choice of terminology was based on a simple majority
rule: we adopted the denomination, ``Feynman diagrams'', of 
the largest community, that of theoretical
physics, which uses the corresponding concept.
Secondly, the historical roots of this notion, especially
in the context of invariant theory, go much further back
in time than Feynman's work.
There is a good account of the history of the diagrammatic notation
in physics and group theory in Chapter 4 of~\cite{Cvitanovic}
to which we refer the reader. This needs, however, to be complemented by the
following pieces of information.

Feynman diagrams, as known to physicists,
seem to have first appeared {\em in print} in~\cite{Dyson},
with due credit to the previously unpublished work of R.~P.~Feynman.
However, the idea of using discrete combinatorial structures 
(e.g.~graphs) in order to describe the outcome of repeated
applications of differential operators with polynomial coefficients
(called ``operandators'' most probably by Sylvester)
such as the polarization, the Omega, and the Aronhold processes of invariant
theory goes back to A.~Cayley~\cite{Cayley}.
Our diagrammatic approach is a presentation, guided by
modern physical wisdom, of the original work of
Sylvester~\cite{Sylvester} and Clifford~\cite{Clifford} (see
also~\cite{Kempe}). It is remarkable that Clifford used what would
now be called Fermionic or Berezin integration
to explain the translation from graphs to actual covariants.
The diagrams we used are a direct visualization of the classical
symbolic notation: arrows correspond to bracket factors,
and each vertex corresponds to a symbolic letter to be repeated
a number of times equal to the degree of the vertex.
There is however one extremely powerful extra feature of the 19th century
methods, in comparison to the physicists' ``diagrammar''.
It is the realization, by Aronhold and Clebsch~\cite{Aronhold, Clebsch},
that one can do all the calculations while pretending that the
ground forms are {\em powers of linear forms}.
This is {\em the} main obstacle lying before a modern who would like
to understand the classics.
Nevertheless, this obstacle can be easily overcome, for instance
by using the umbral methods developed by Rota and his school~\cite{Kung}.
Here is another way to rigorously justify this simplification, which
we believe offers more flexibility, like for instance the possibility
of iteration, that of mixed interpretation of some variables
as true and others as symbolic within the same computation,
as well as that of using, in intermediate steps,
the same symbolic letter a number of times
exceding the degree of the form.
It goes as follows:
do all the needed calculations {\em as stated} with symbolic letters
$a$, $b$, $c$\ldots
considered as formal indeterminates,
and afterwards, act on the resulting expression with
the product 
\[ Q(\frac{\partial}{\partial a}) \, 
Q(\frac{\partial}{\partial b}) \, 
Q(\frac{\partial}{\partial c})\ldots \] 
where $Q$ is the ground form under consideration.
This will insert the coefficients of the form
at the right place in the symbolical expression, i.e., at the right
vertex of the diagram.

Diagrams, in the context of classical invariant theory,
were reintroduced in the work of Olver and Shakiban~\cite{OShakiban}
which is a slightly different formalism because of a normal ordering
procedure explained in Chapter 6 of~\cite{Olver}.
Finally, although somewhat atypical, the interesting 
pedagogical work of computer graphics pioneer J.~F.~Blinn~\cite{Blinn},
who was inspired by the book~\cite{Stedman}, deserves to be mentioned.

\medskip 

\parbox{12.5cm}{\small 
{\sc Acknowledgements:}
The first author would like to thank David C.~Brydges and 
Joel S.~Feldman for their kind invitation to visit the University of British Columbia. 
The second author would like to thank Prof.~James Carrell and 
and the University of British Columbia for financial assistance. 
We are indebted to the authors of the package Macaulay-2.}

\bibliographystyle{plain}
\bibliography{../../BIBLIO/ref1,../../BIBLIO/ref2,../../BIBLIO/ref3}

\vskip 2cm
\parbox{12cm}{\small 
{\sc Abdelmalek Abdesselam}
\begin{flushleft}
Department of Mathematics\\
University of British Columbia,\\
1984 Mathematics Road,\\
Vancouver, B.~C.~, V6T 1Z2, CANADA.
\end{flushleft}

\noindent 
and

\begin{flushleft}
LAGA, Institut Galil\'ee, CNRS UMR 7539,\\
Universit{\'e} Paris XIII,\\
99 Avenue J.B. Cl{\'e}ment\\
F93430 Villetaneuse, FRANCE.
\end{flushleft}

\bigskip
abdessel@math.ubc.ca}

\vskip 2cm

\parbox{12cm}{\small 
{\sc Jaydeep Chipalkatti}
\begin{flushleft}
Department of Mathematics\\
University of British Columbia,\\
1984 Mathematics Road,\\
Vancouver, B.C., V6T 1Z2, CANADA.
\end{flushleft}

\bigskip
jaydeep@math.ubc.ca}

\end{document}